\documentclass[11pt]{article}
\usepackage{amsmath}
\usepackage{amssymb}
\usepackage{textcomp}
\usepackage[document]{ragged2e}
\usepackage{setspace}
\usepackage{comment}
\usepackage{slashed}
\usepackage{etoolbox}
\usepackage[affil-it]{authblk}

\usepackage{tikz}
\usetikzlibrary{positioning}
\usepackage{dynkin-diagrams}
\usepackage[pdfencoding=auto, psdextra]{hyperref}
 \hypersetup{
     colorlinks=true,
     linkcolor=blue,
     citecolor = blue,      
     urlcolor=blue,
     }
\newenvironment{changemargin}[2]{
\begin{list}{}{
\setlength{\topsep}{0pt}
\setlength{\leftmargin}{#1}
\setlength{\rightmargin}{#2}
\setlength{\listparindent}{\parindent}
\setlength{\itemindent}{\parindent}
\setlength{\parsep}{\parskip}
}
\item[]}{\end{list}}

\newcommand{\bb}[1]{\mathbb{#1}}
\newcommand{\alg}[1]{\mathfrak{#1}}
\newcommand{\spin}{\mathfrak{spin}}
\newcommand{\su}{\mathfrak{su}}
\newcommand{\so}{\mathfrak{so}}
\newcommand{\grp}[1]{\operatorname{#1}}

\title{Exploring Triality Explicitly: \newline Convenient bases for $\grp{SO}(8)$, $\grp{Spin}(1,7)$, and $ \grp{G}_2$}
\author{Craig M$^{\mathrm{c}}$Rae\thanks{Electronic address: \texttt{mcraec3@myumanitoba.ca}}}
\affil{University of Manitoba \\ Winnipeg, MB}

\begin{document}

\maketitle
  \pagenumbering{gobble}
%  \newpage
  \pagenumbering{arabic} 
\baselineskip=1.2\baselineskip
\setlength{\parindent}{1.5em}
\setlength{\parskip}{0.5ex}

\begin{changemargin}{-1cm}{-1cm}

\begin{abstract}
\normalsize The property of triality only appears in one linear simple Lie algebra: $D_4$, a.k.a. $\mathfrak{so}(8, \mathbb{C})$. In this paper several convenient representations and bases are constructed in order to facilitate the exploration of the three fold symmetry known as the triality of representations. In particular the three $8$ dimensional representations for the Euclidean and Lorentzian real forms of $\mathfrak{so}(8,\mathbb{C})$ are constructed, and the maps between representations are given in each case, respectively. Standard results are presented in a modern language, such as how $\mathfrak{g}_2$ arises as the intersection of non-conjugate $\mathfrak{spin}(7,\mathbb{R})$ sub-algebras, and as the stabilizer of the outer automorphism group of $\mathfrak{so}(8)$. 

It is argued that $\mathfrak{spin}(1,7)$ is in some sense the more natural stage for triality to play out upon, as triality can be made to act diagonally in this case, leading to a $\mathbb{Z}_3$ grading of the Lie algebra and subsequently a Jordan algebra. Once this is understood a short discussion is had about obstacles to a theory of triality which attempt to explain the three generations of matter via some form of triality.
\end{abstract}
\end{changemargin}
\begin{changemargin}{-2cm}{-2cm}

\newpage 
\tableofcontents

\newpage
\section{The Triality Automorphism}
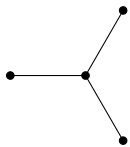
\begin{figure}[!h]
\centering{
\begin{tikzpicture}
\dynkin[edge length=1cm]D4
\end{tikzpicture}
\caption{The Dynkin diagram $D_4$ of $\alg{so}(8, \mathbb{C})$}
\label{fig:D4}
}
\end{figure}

The Dynkin diagram for the the Lie algebra $\alg{so}(8)$, also known as $D_4$ in the Dynkin classification, is the most symmetric diagram for any complex Lie algebra with finite dimensional representations. The diagram shares the symmetries of an equilateral triangle in the plane, the symmetric group $S_3$. Since symmetries of an algebra's Dynkin diagram correspond to outer automorphisms of its representations \cite{FultonHarris}, we find irreducible representations of $\alg{so}(8)$ all fall in to one of two kinds: either they are a singlet under the outer automorphism group, like the adjoint representation, thought of as the middle point of the diagram, or they are one of a triplet of representations enjoying a ménage à trois, mapped to one another via the outer automorphism group of the algebra. Often the discussion of this fact is left here, as a mathematical curiosity, and explored only in greater abstraction. While this is insightful,\footnote{Such as in this wonderful paper \cite{VARADARAJAN2001}, wherein I first learned many of the stated facts on display here.} it can leave a curious mathematician wondering how they might actually see the triality at play: what does the performance of these automorphisms actually \textit{look like}? In what follows an explicit mapping between the three $\mathbf{8}$ dimensional representations of $\alg{so}(8, \mathbb{R})$ is given, in order for the reader to have a more concrete understanding of this symmetry between the representations. 
\section{Compact Real Form: \texorpdfstring{$\so(8,\mathbb{R})$}{so(8,R)}}
\subsection{Vector Representation}
The defining representation of $\grp{SO}(8, \mathbb{R})$, is the $8 \times 8$ real orthogonal matrices $R$ satisfying:
\begin{equation}
    R^\intercal R = \mathbb{I}_8.
\end{equation}
It is easy to show using the standard techniques that the Lie algebra $\mathfrak{so}(8, \mathbb{R})$ of this representation is identified with the $28$ dimensional space of real anti-symmetric $8 \times 8$ matrices. A convenient basis for the algebra is given by matrices $V_{ij}$, with $0 \leq i < j \leq 7 $, with a $1$ in the $i^\mathrm{th}$ row and $j^\mathrm{th}$ column, anti-symmetric, and zero everywhere else. This is a simple way to generate the $28$ planes of rotation in a Cartesian coordinate system for $8$ dimensional space. A condensed formula for the components of any generator of rotations in this basis is
\begin{equation}
\mathbf{8}_V = \left(V_{ij}\right)_{a,b} = \delta^{i}_{a}\delta^{j}_{b}-\delta^{j}_{a}\delta^{i}_{b}
\end{equation}
with $\delta$ the Kronecker delta. Here $V$ stands for the \textit{vector} representation, as we will see there are also two spinor representations\footnote{For brevity a particular technicality is avoided by staying at the level of the algebra. The group representations one finds by exponentiation of the spin algebras built in the following section are technically speaking not \textit{spin} reps. They are what are called \textit{semi-spin} representations. As a double cover, a spin group should have a $2:1$ homomorphism to the corresponding special orthogonal group. In particular the center must be double covered. For $SO(2n)$ the center is $\mathbb{Z}_2$, and the double covers of this group are the Klein four group $\mathbb{Z}_2 \times \mathbb{Z}_2$, and the cyclic group of order four $C_4$. When $2n \grp{mod} 4 = 2$ the center is $C_4$, while for $2n \grp{mod} 4 = 0$ the center is the Klein group. The latter case causes an obstruction for faithful spin irreps, because irreps of Abelian groups are faithful if and only if the group is cyclic \cite{Huppert1998}. What this means is that while the group $\grp{Spin}(8)$ in principle has a center with four elements, any irrep we find will be missing `half' of the center. One may think of this as modding out by the `other' $\mathbb{Z}_2$ instead of the usual one done to acquire the orthogonal group. This means the only true spin reps are reducible reps such as $L\oplus R$ where the volume element $\omega$ is central.} which act on an $8$ dimensional real space as well. Finally we will alter this simple basis of rotation generators as follows: flip $V_{1,5} \mapsto -V_{1,5}$ and $V_{2,6} \mapsto -V_{2,6}$, all else is left the same  (this should be no more mysterious than choosing a right hand rule in three dimensions).

\subsection{Real Chiral Spinor Representations}
To build bases of the spinor representations of $\operatorname{Spin}(8)$ we can start with a more familiar Clifford algebra and build up to $C\ell(8)$. First with the usual gamma matrices for $\alg{spin}(1,3)$ in the Dirac basis,\footnote{These can be found in any quantum field theory text, see for example \cite{Schwartz2014}.} we may build a completely imaginary\footnote{Should one wish, one finds a totally real Clifford algebra by dividing each of these by an $i$, at the cost of making the default signature totally negative: $C\ell(0,7)$.} matrix basis $g_i$ for $C\ell(7)$. 
\begin{equation}
\label{eq:so(7)}
\begin{split}
g_1 &= i\sigma_z \otimes \gamma_1 \gamma_3 = i \begin{pmatrix}
\gamma_1 \gamma_3 & 0 \\
0 & -\gamma_1 \gamma_3
\end{pmatrix},\\
g_2 &= i\sigma_z \otimes \gamma_3, \quad g_3 = i \sigma_z \otimes \gamma_1, \\
g_4 &= -\sigma_y \otimes \mathbb{I}_4, \quad g_5 = \sigma_x \otimes \gamma_5 \gamma_2, \\
g_6 &= i\sigma_x \otimes \gamma_0 \gamma_5, \quad g_7 = \sigma_x \otimes \gamma_2 \gamma_0.
\end{split}
\end{equation}
With this basis, we can construct a new completely real basis $\Gamma_\mu$ for $C\ell(8,0)$
\begin{equation}
\Gamma_0 = \begin{pmatrix}
0 & \mathbb{I}_8 \\
\mathbb{I}_8 & 0
\end{pmatrix}, \quad \Gamma_{\mu} = -i \begin{pmatrix}
0 & g_\mu \\
-g_\mu & 0
\end{pmatrix} \quad 1 \leq \mu \leq 7.
\end{equation}
These satisfy the standard Euclidean Clifford algebra $\{\Gamma_i, \Gamma_j\} = 2\delta_{ij}$ for $0 \leq i,j \leq 7$, where $\{,\}$ is the anti-commutator. As usual,\footnote{For a great background and introduction to the relationship between Clifford algebras and spin groups see \cite{ClifNSpinor}. } the $\mathfrak{spin}(8)$ algebra can be generated by second degree elements of the Clifford algebra, and conveniently from the $\Gamma_i$ above we can see the generators in this basis will all be block diagonal
\begin{equation}
\label{LRdef}
\mathbf{8}_L \oplus \mathbf{8}_R = L_{ij} \oplus R^\prime_{ij}=\operatorname{span}\{\Gamma_i \Gamma_j/2\}, \quad 0 \leq i < j \leq 7.
\end{equation}
By this notation I mean the blocks of this 28 dimensional basis of our $\mathfrak{spin}(8)$ Lie algebra decompose into: the upper $8 \times 8$ block $L_{ij}$ is the left handed spin representation, and the lower block is the right handed spin representation $R^{\prime}_{ij}$ (the prime is here because we will soon alter the right handed basis to play nicer with triality). Each obey identical commutation relations (the structure functions of the bases match). In this basis they are related via a reflection along the $0^{\mathrm{th}}$ axis, i.e. if one sends $\Gamma_0 \mapsto -\Gamma_0$, then $L \leftrightarrow R^\prime$. This implies that when $i,j > 0$ then $L_{ij} = R^\prime_{ij}$ --- indeed this must be true up to a change of basis as restricting to a seven dimensional subspace should give generators of $\mathfrak{spin}(7)$, for which there is no handedness. The factor of $\frac{1}{2}$ is a matter of convention but still important to keep things consistent with the previously written down vector representation. One can also define a (real) volume element $\omega$ which is the product of the $8$ gamma matrices, anti-commutes with all $8$ gamma matrices, squares to the identity, and can be used to define left and right projection operators. For yet unseen reasons of convenience, we will employ the following change of basis to all the right handed generators:
\begin{equation}
\label{Rredef}
    R_{ij} \> = \>P R^{\prime}_{ij} P^\intercal, \qquad P = \grp{diag}\{-1,1,1,1,1,1,1,1\}.
\end{equation}
Recall the indices are referring to a particular matrix generator, not the elements of a matrix. It will also be of note later that the change of basis has $\det = -1$, and so is not a member of $\grp{Spin}(8)$. Finally, as above, we will use the basis for $\mathbf{8}_L $ and $\mathbf{8}_R$ given in Eqs.~(\ref{LRdef}, \ref{Rredef}), but with the same alteration as the vector case where we flip the default sign of $L_{1,5},\>R_{1,5},\>L_{2,6}$, and $R_{2,6}$.

\subsection{Triality Explicitly}
That there are three unique $8$ dimensional representations of $\mathfrak{spin}(8)$ is noteworthy in its own right, but triality is more interesting than this. Playing around with these three ways of performing rotations on an $8$ dimensional vector space one may notice each is a $28$-dimensional vector space of real, antisymmetric matrices. But the vector space of $8\times8$ anti-symmetric matrices is precisely $28$ dimensional --- the generators of the $\mathbf{8}_V$, $\mathbf{8}_L$, and $\mathbf{8}_R$ representations of $\mathfrak{spin}(8)$ all span \textit{precisely} the same space. As such there must be a way to identify them with one another, i.e it should be true any $L_{ij}$ can be written in terms of some linear combination of $V_{ij}$'s, etc. What is this map? 

\subsubsection{The Triality Map}
For the representations to be distinct, the matrices cannot be merely a change of basis away from one another, which one can verify explicitly is not the case. This means the automorphism must act foremost at the level of the $28$ dimensional vector space of the Lie algebra. In the bases given, the triality map is as follows: collect the $28$ basis elements into four ordered sets of seven. In particular define the sets (thought of as seven quartets being simultaneously acted on):\footnote{Note no element of any row commutes with another in its row, while all elements of a column commute.}
\begin{equation}
\label{TrialityQuartets}
\begin{split}
\vec{a}_V &= \{V_{01},V_{02},V_{03},V_{04},V_{05},V_{06},V_{07}\}, \\
\vec{b}_V &= \{V_{23},V_{57},V_{12},V_{37},V_{36},V_{17},V_{25}\},\\
\vec{c}_V &= \{V_{45},V_{13},V_{47},V_{15},V_{14},V_{24},V_{16}\},\\
\vec{d}_V &= \{V_{67},V_{46},V_{56},V_{26},V_{27},V_{35},V_{34}\},
\end{split}
\end{equation}
and analogously for $L$ and $R$. The following map on this $28$ dimensional space is a Lie algebra homomorphism which acts cyclically on the $V, L, R$ representations:
\begin{equation}
\begin{pmatrix}
\vec{a} \\
\vec{b} \\
\vec{c} \\
\vec{d}
\end{pmatrix}_{\hspace{-0.5em} L} = \frac{1}{2}\begin{pmatrix}
-1 & -1 & 1 & 1 \\
1 & 1 & 1 & 1 \\
-1 & 1 & 1 & -1 \\
-1 & 1 & -1 & 1 
\end{pmatrix}\begin{pmatrix}
\vec{a} \\
\vec{b} \\
\vec{c} \\
\vec{d}
\end{pmatrix}_{\hspace{-0.5em} V}.
\end{equation}
Let us call the matrix given (including the factor of $\frac{1}{2}$) $H$. If one applies this map to all seven sets of four generators in the vector representation, they will find they recover precisely the generators of the left handed spinor representation. Applying it again sends the left handed spinor rep to the right handed spinor rep defined above, and applying it a third time returns to the vector rep. This is to say $H^3 = \bb{I}_4$ which of course implies $H^2 = H^{-1}$. It is worth pointing out one can unpack $H$ if they so choose to be a $28$ dimensional square matrix, this is merely a compact way of writing this down. Let us also define another matrix $K$:
\begin{equation}
K = \begin{pmatrix}
-1 & 0 & 0 & 0 \\
0 & 1 & 0 & 0 \\
0 & 0 & 1 & 0 \\
0 & 0 & 0 & 1 
\end{pmatrix}.
\end{equation}
This matrix acting on the same set of generators as $H$, simply flips the $\vec{a}$ entries in all the quartets of bases. As noted above this may be understood to be a reflection through the $0^{th}$ axis. This action on $\mathbf{8}_L$ or $\mathbf{8}_R$ will send one to the other (though in a slightly less satisfying way because of the change of basis introduced for our $\mathbf{8}_R$ basis\footnote{Were we to remove this change of basis it will instead show up as needed in the action of $H$, so it is left a problem of $K$ where it is easiest to keep track of.}). The action of $K$ on the vector rep leaves us in place, at least at the level of the representation, however it will still have the effect of essentially performing a change of basis reflection along the $0^{th}$ axis here too. 

We have two operations on the Lie algebra of $\spin(8,\mathbb{R})$ which are algebra homomorphisms, and play relatively nicely with our chosen bases. The action of $K$, with $K^2 = 1$, and the action of $H$, where $H^3 = 1$. At the level of representations $H$ is a rotation of the Dynkin diagram by $2 \pi /3$ and $K$ is a reflection which fixes the vector node. The only complication being that while $K$ certainly maps between representations, there is no nice way to have $H$ and $K$ simultaneously map exactly between the same bases, and so applications of $K$ must be followed up by a change of basis $P$ to return to the bases defined above. Inspection of how $K$ and $H$ interact reveal unexpectedly that $K$ and $H$ generate a representation of the permutation group $S_3$, where the thing being permuted of course is which $8$-dimensional representation of $\mathfrak{spin}(8)$ we are in.

Recall that this is a real symmetry of all three representations, even though this basis emphasizes the usual duality in the spinor representations. The fact that the algebra outer automorphisms form the group $S_3$ means there is just as good a permutation swapping $V$ and $L$ while leaving $R$ alone --- contained in the triality are three dualities. 

\subsection{Non-Conjugate of \texorpdfstring{$\alg{spin}$(7) Sub-Algebra's of $\alg{spin}$(8)}{spin(7) sub-algebra's of spin(8)}}
In the representations given, there are many ways one can restrict the $28$ generators to some $\alg{spin}(7)$ sub-algebra. Generally, removing rotations involving any one axis will have this effect, and thinking of the vector rep for example, one does not really consider different choices of axes being removed as `different' $\alg{spin}(7)$ sub-algebras, since one could simply rotate by some $\grp{SO}(8)$ transform and change which chosen axis is being dropped. The question then arises if there exist `distinct' $\alg{spin}(7)$ sub-algebra's of $\alg{spin}(8)$, which cannot be merely mapped to one another under the action of $\grp{Spin}(8)$. To a mathematician this is the question of how many conjugacy classes of $\alg{spin}(7)$ there are inside of $\alg{spin}(8)$. The answer is that there are precisely $3$ conjugacy classes of $\alg{spin}(7)$ inside of $\alg{spin}(8)$, and this is related precisely to the triality: each of these conjugacy classes are mapped to one another under the action of the triality automorphism.\footnote{Relating to footnote $4$, Varadarjan \cite{VARADARAJAN2001} has a wonderful discussion of how one can understand the three non-conjugate $\grp{Spin}(7)$ subgroups as essentially being non-conjugate because they each inherit a different non-trivial element of the center from the full $\grp{Spin}(8)$, and then going on to discuss how triality acts on the center of the group via the outer automorphism group of the Klein group, which is also $S_3$.} In the course of investigating this, instead of working in one representation and finding the non-conjugate $\alg{spin}(7)$ sub-algebra's thereof, it is easier to simultaneously hold our three $\alg{spin}(8)$ representations in mind, and remove the same axis in each of them, rather than think of them as three sub-algebra's of one representation. The math of course is the same in either case.

With the above expounded the selection of a basis for any three non-conjugate $\alg{spin}(7)$ sub-algebras is as simple as picking an axis to remove (an index to disallow) from our sets of rotation generators related by the action of $H$: $V_{ij}$, $L_{ij}$, and $R_{ij}$. Once this choice is made one will have three representations of $\alg{spin}(7)$: a vector representation and two spinor representations. As $\alg{spin}(7)$ has only a single spinor representation these two sets of generators \textit{must} be linearly related, however we know they are related by a reflection $P$, and so no $\grp{Spin}(8)$ transform could have mapped these sub-algebras to one another. It is somewhat more obvious to see that there is no way to map the restricted set of vector generators to either of the spinor cases, since the algebras need no longer span an identical subspace of anti-symmetric matrices, and it is well known that the spinor and vector representations of $\spin(7)$ are distinct irreps. We'll construct and utilize these restrictions in the next section.

\subsection{The Intersection of Non-Conjugate \texorpdfstring{$\alg{spin}(7)$ Sub-Algebra's: $\alg{g}_2$}{spin(7) sub-algebra's: g2}}
The triality of $\alg{spin}(8)$ gives rise to three non-conjugate $\alg{spin}(7)$ sub-algebras. Of course one cannot fit three $\alg{spin}(7)$ sub-algebras neatly inside $\alg{spin}(8)$, and so a natural question is how much these three conjugacy classes `overlap', that is to say what is the subspace of $\alg{spin}(8)$ which is `fixed' under the outer automorphism group. This is given by the intersection of the distinct conjugacy classes. The answer to this is that the subspace $\alg{spin}(8)$ which is simultaneously in the span of all three $\alg{spin}(7)$ sub-algebra's, is $14$ dimensional. This 14 dimensional sub-space is mapped to itself under triality, which means it must be closed under the Lie bracket.\footnote{For a less complicated analogous situation, consider the two conjugacy classes of $\alg{su}(2)$ in $\alg{su}(3)$: the vector representation spanned by $\lambda_5,\lambda_7,\lambda_2$ and the spinor representation spanned by $\lambda_1, \lambda_2, \lambda_3$. Their intersection is an $\so(2)$.} Which semi-simple Lie group is this, the subgroup of $\grp{Spin}(8)$ serving as the stabilizer of the outer automorphism group? The answer is the exceptional Lie group $G_2$! Let us work this through in our bases. 

\subsubsection{Finding \texorpdfstring{$\alg{g}_2$ in $\spin(8)$}{g2 in spin(8)}}
For our case let us remove the $0^{th}$ axis, since we understand its effects particularly well. This restricts all three of our sets of generators to have $ij > 0$, and so we find three sets of $21$ dimensional $\spin(7)$ algebras. As noted at the start, from the construction of the Clifford algebra we know that the restrictions of $L$ and $R$ to this particular choice of sub-algebra are related by application of the change of basis $P$. This may seem ad hoc but one must remember the $P$ is there to play nice with triality, not merely to make the sub-algebras non-conjugate. As noted $P$ is not a member of $\grp{Spin}(8)$ and so these sub-algebras are non-conjugate. The sub-algebra coming from the vector representation is also not conjugate to either of these, which is easy to see as discussed above. 

With these restricted bases from $L$, $R$, and $V$, next we need to find the intersection of these sub-algebras. This is merely the question of, within the vector space of $\spin(8)$, what is the set of vectors which are in the span of all three of these $\spin(7)$ subspaces?  This is an algebra problem: let a generic vector in one basis equal a generic vector in another, giving constraints on the coefficients. In particular we demand
\begin{equation}
    \sum_{0<i<j\leq 7} a_{ij} V_{ij} = \sum_{0<i<j\leq 7} b_{ij} L_{ij}.
\end{equation}
This imposes the following restrictions on the coefficients:
\begin{equation}
\label{interseven}
\begin{split}
    b_{1,2} = b_{4,7} + b_{5,6}, \quad &b_{1,3} = -b_{4,6} + b_{5,7}, \\
    b_{1,4} = -b_{2,7} + b_{3,6}, \quad &b_{1,5} = -b_{2,6} + b_{3,7}, \\
    b_{1,6} = b_{2,5} - b_{3,4}, \quad &b_{1,7} = b_{2,4} + b_{3,5}, \\
    b_{2,3} = b_{4,5} + b_{6,7}, \quad &a_{i,j} = b_{i,j} \quad \forall i,j.\\
\end{split}
\end{equation}
For those counting along at home, this was $42$ parameters, with $7+21$ constraints, leaving us with $14$ free dimensions. If one takes the remaining $R$ basis, and takes an arbitrary $21$ dimensional vector $\sum_{ij} b_{ij} R_{ij}$, and imposes the restrictions from Eqs.~(\ref{interseven}), the will find they are restricted to exactly the same subspace, and so can verify demanding any two of the intersection is as good as demanding all three. We have found our $14$ dimensional subspace fixed by the outer automorphism group! In the bases given, it is the following sub-algebra:
\begin{equation*}
\hspace{-0.5em}
\begin{pmatrix}
 0 & 0 & 0 & 0 & 0 & 0 & 0 & 0 \\
 0 & 0 & b_{4,7}+b_{5,6} & b_{5,7}-b_{4,6} & b_{3,6}-b_{2,7} & b_{2,6}-b_{3,7} & b_{2,5}-b_{3,4} & b_{2,4}+b_{3,5} \\
 0 & -b_{4,7}-b_{5,6} & 0 & b_{4,5}+b_{6,7} & b_{2,4} & b_{2,5} & -b_{2,6} & b_{2,7} \\
 0 & b_{4,6}-b_{5,7} & -b_{4,5}-b_{6,7} & 0 & b_{3,4} & b_{3,5} & b_{3,6} & b_{3,7} \\
 0 & b_{2,7}-b_{3,6} & -b_{2,4} & -b_{3,4} & 0 & b_{4,5} & b_{4,6} & b_{4,7} \\
 0 & b_{3,7}-b_{2,6} & -b_{2,5} & -b_{3,5} & -b_{4,5} & 0 & b_{5,6} & b_{5,7} \\
 0 & b_{3,4}-b_{2,5} & b_{2,6} & -b_{3,6} & -b_{4,6} & -b_{5,6} & 0 & b_{6,7} \\
 0 & -b_{2,4}-b_{3,5} & -b_{2,7} & -b_{3,7} & -b_{4,7} & -b_{5,7} & -b_{6,7} & 0 
 \end{pmatrix}
\end{equation*}
Notice that since we chose to remove the $0^{th}$ axis for this, and we must intersect the vector representation, entries must be null in that row and column. All things considered treating these coefficients as $0$'s and $1$'s would serve for a relatively nice basis, however it would satisfy neither normality or orthogonality. A more convenient orthonormal basis $\Lambda_i$ for this sub-algebra may be given instead by (omitting the $0^{th}$ row and column):

\begin{equation}
\begin{split}
 \sum_{i=1}^{14}\theta_i \Lambda_i &= 
\frac{1}{2}\begin{pmatrix}
 0 & 0 & 0 & 0 & 0 & 0 & 0 \\
 0 & 0 & 0 & \theta_6 & \theta_7 & -\theta_5 & \theta_4 \\
 0 & 0 & 0 & \theta_7 & -\theta_6 & \theta_4 & \theta_5 \\
 0 & -\theta_6 & -\theta_7 & 0 & \theta_3 & \theta_1 & \theta_2 \\
 0 & -\theta_7 & \theta_6 & -\theta_3 & 0 & -\theta_2 & \theta_1 \\
 0 & \theta_5 & -\theta_4 & -\theta_1 & \theta_2 & 0 & -\theta_3 \\
 0 & -\theta_4 & -\theta_5 & -\theta_2 & -\theta_1 & \theta_3 & 0 \\
\end{pmatrix} \\
& +\frac{1}{2\sqrt{3}}
\begin{pmatrix}
 0 & -2 \theta_9 & -2 \theta_{10} & -2 \theta_{11} & 2 \theta_{12} & -2 \theta_{14} & -2 \theta_{13} \\
 2 \theta_9 & 0 & -2 \theta_{8} & -\theta_{13} & -\theta_{14} & -\theta_{12} & \theta_{11} \\
 2 \theta_{10} & 2 \theta_{8} & 0 & \theta_{14} & -\theta_{13} & -\theta_{11} & -\theta_{12} \\
 2 \theta_{11} & \theta_{13} & -\theta_{14} & 0 & -\theta_{8} & \theta_{10} & -\theta_9 \\
 -2 \theta_{12} & \theta_{14} & \theta_{13} & \theta_{8} & 0 & -\theta_9 & -\theta_{10} \\
 2 \theta_{14} & \theta_{12} & \theta_{11} & -\theta_{10} & \theta_9 & 0 & -\theta_{8} \\
 2 \theta_{13} & -\theta_{11} & \theta_{12} & \theta_9 & \theta_{10} & \theta_{8} & 0 \\
\end{pmatrix}.
\end{split}
\end{equation}
Here $\theta_i$ are $14$ real variables (coefficients). What is appealing about this basis is that $\Lambda_1$ through $\Lambda_8$ form precisely a standard $\mathfrak{su}(3)$ sub-algebra, so this limit is simply taken by ignoring the last $6$ generators. The basis is also nice in two other aspects. Firstly that it splits into 7 generators which are `like' the first seven Gell-Mann matrices and 7 generators which are `like' the eighth Gell-Mann matrix $\lambda_8$.\footnote{The word `like' here should be read as `their spectra are analogous'.} Second, if one swaps the names of $\Lambda_{8}$ and $\Lambda_{10}$, making the $\su(3)$ sub-algebra slightly less obvious, then its true that $[\Lambda_i,\Lambda_{i+7}] = 0$. One can verify in this orthonormal basis that the Lie bracket is closed, and since $G_2$ is the only simple fourteen dimensional Lie algebra, we must have found an instance of (the compact real form of) $\alg{g}_2$! 

\subsubsection{Finding \texorpdfstring{$\su(3) $ in $\alg{g}_2$}{su(3) in g2}}
Descending even further to obtain the usual representation's of $SU(3)$, we can simply apply the following special unitary transform to the $8$ basis elements which form the $\mathfrak{su}(3)$ sub-algebra:

\begin{equation}
U = \frac{1}{\sqrt{2}}\begin{pmatrix}
 \sqrt{2} & 0 & 0 & 0 & 0 & 0 & 0 \\
 0 & 0 & 0 & 0 & 0 & 1 & -i \\
 0 & 0 & 0 & -i & -1 & 0 & 0 \\
 0 & -1 & -i & 0 & 0 & 0 & 0 \\
 0 & 0 & 0 & 0 & 0 & -i & 1 \\
 0 & 0 & 0 & 1 & i & 0 & 0 \\
 0 & i & 1 & 0 & 0 & 0 & 0 \\
\end{pmatrix}.
\end{equation}
After application of this unitary, the real $\mathfrak{su}(3)$ sub-algebra of $\alg{g}_2$ decomposes as:
\begin{equation}
    U \Lambda_i U^\dagger = \begin{pmatrix}
        0 & 0 & 0 \\
        0 & \lambda_i & 0 \\
        0 & 0 & -\lambda_i^\intercal
    \end{pmatrix},
\end{equation}
which is precisely $ \mathbf{1} \oplus \mathbf{3} \oplus \overline{\mathbf{3}} $. It should be noted the physics convention is gained by slapping an $i$ on all $\Lambda_i$.

\section{Diagonalizing Triality and \texorpdfstring{$\so(8, \mathbb{C})$}{so(8,C)}}
An easier though less geometric way of finding $G_2$ through its relationship to triality, is inspecting the eigenvectors of the triality matrix $H$: 
\begin{equation}
\label{Heigens}
    \begin{split}
    &\left| 1\right\rangle_3 = \frac{1}{\sqrt{2}}\begin{pmatrix}0 \\ 1 \\ 0 \\ 1\end{pmatrix}, \quad \left| 1\right\rangle_8 = \frac{1}{\sqrt{6}}\begin{pmatrix} 0 \\1 \\2 \\ -1\end{pmatrix}, \quad \left|e^{i 2 \pi /3}\right\rangle = \frac{1}{\sqrt{6}}\begin{pmatrix} i \sqrt{3} \\ -1 \\ 1 \\ 1 \end{pmatrix} = \left| e^{-i 2 \pi /3}\right\rangle^*
    \end{split}
\end{equation}
The vectors are labeled by their eigenvalues under triality, and $^*$ is complex conjugation. The unitary change of basis matrix defined by 
\begin{equation}
    U = \left(\left| 1\right\rangle_3 \> \left| 1\right\rangle_8 \> \left| e^{i 2 \pi /3}\right\rangle \> \left| e^{-i 2 \pi /3}\right\rangle\right),
\end{equation}
diagonalizes $H = U D U^\dagger$ where $D = \operatorname{diag}\{1,1,e^{i 2 \pi /3},e^{-i 2 \pi /3}\}$. Recall $H$, and hence $U$, act on quartets of $\so(8, \mathbb{R})$ generators, so this change of basis gives us complex combinations of our generators. This means we generically can only anticipate we find ourselves with a basis for $\so(8,\mathbb{C})$. However once we have done this change of basis our 28 generators are then conveniently arranged into a $14$ dimensional sub-algebra which have a trivial orbit under the action of $\mathrm{Out}(\so(8, \mathbb{C}))$, and two sets of seven generators, which accrue the given eigenvalues of conjugate third roots of unity under the operation of triality. Let us refer to the seven dimensional subspace with eigenvalue $e^{i 2 \pi /3}$ the `right handed' subspace, and the seven dimensional subspace with eigenvalue $e^{-i 2 \pi /3}$ the `left handed' subspace.

The fourteen dimensional outer automorphism invariant sub-algebra is of course a basis for $\alg{g}_2$: not only does triality map the $\alg{g}_2$ sub-algebra to itself, from the eigenvalues of $H$ we can see it fixes each of those elements in place. Not only is the outer automorphism orbit of the $\alg{g}_2$ sub-algebra closed, it is in fact trivial. This is also the origin of the $3,8$ subscripts on the eigenvectors. The eigenvalue $1$ shows up with multiplicity of two and so the $1$-eigenspace can be decomposed orthogonally into $7$ bases which `look like' $\lambda_3$ and $7$ which `look like' $\lambda_8$, which should be familiar from earlier.\footnote{It is worth noting though there is still an $\su(3)$ sub-algebra, it is not as obvious here. This is due to the fact that for any quartet upon which triality acts in the original basis, the orthogonal eigen-decomposition of the $1$ eigenspace can in principle be different for each quartet, and would need to be done here for the standard $\su(3)$ commutation relations to arise. For simplicity this is not done here.} The remaining $14$ generators which are non-trivially affected by triality do not form a closed sub-algebra. Since triality is an algebra homomorphism, this is easily seen by inspecting the action of triality on the commutator: the commutator of any two vectors in the right handed eigenspace can be seen to be in the left handed eigenspace and vice versa. Furthermore the commutator of two vectors, one of each handed eigenspace, must have eigenvalue $1$ under triality and so belong inside the $\alg{g}_2$ sub-algebra. 

The given basis lifts the veil on triality perhaps as much as a good basis can, but this does not refute its intrigue. We find here there are $28$ generators of $\so(8,\mathbb{C})$ and they may be indexed by which third root of unity they are multiplied by under the triality automorphism: $\{1, e^{i 2 \pi /3},e^{-i 2 \pi /3}\}$ where the eigenspaces are $14, 7,$ and $ 7$ dimensional respectively. Another way to say this is that the distinction between the three $\mathbf{8}$ dimensional representations is merely a particular phase applied to those two sets of 7 generators, making the automorphism far more similar to a negation (second root of unity) being applied to a subset of generators, which is the structure of most Lie algebra dualities. A final interesting thing to inspect is the action of $K$ in this basis. Applying the change of basis $U$ we find:
\begin{equation}
    K^\prime = U^\dagger K U = \begin{pmatrix} 1 & 0 & 0 & 0 \\
    0 & 1 & 0 & 0 \\
    0 & 0 & 0 & 1 \\
    0 & 0 & 1 & 0 \\
    \end{pmatrix},
\end{equation}
$K^\prime$ unsurprisingly fixes the $\alg{g}_2$ sub-algebra, and appropriately, `swaps' the two seven dimensional eigenspaces, mirroring its behaviour upon the Dynkin diagram.

One might still wish consider this a basis of $\so(8,\mathbb{R})$ sitting somehow skew in its complexification, as real combinations of these bases will still preserve the standard bi-linear inner product given by $g = \mathbb{I}_8$.\footnote{It is noteworthy that in this basis the vector generators only will also preserve the standard Hermitian sesqui-linear product $ h =\mathbb{I}_{1,7}$} However the analogy breaks beyond this as those Lie algebra bases outside the $\alg{g}_2$ sub-algebra (the generators upon which triality acts non-trivially) are null in this basis: their squares are traceless and the Killing form is not diagonal. Of these two sets of seven generators spanning their `handed' eigenspaces respectively, each basis will commute with precisely one from the set with conjugate eigenvalue, and each of these basis elements is orthogonal to all other bases aside from this `conjugate sibling' which it commutes with, and so it seems this diagonalization has naturally put us in a form where the handed bases take the form of raising and lowering operators $L_{xy} \pm iL_{wz} $ etc. As such we are certainly no longer strictly utilizing a basis of $\so(8,\mathbb{R})$. In fact the trace of the normalized Killing form, originally $-28$, here is now $-14$, which seems to imply we are in a basis which may be more naturally considered to be the algebra of $\so(1,7)$. 

\section{Lorentzian case: \texorpdfstring{$\spin(1,7, \mathbb{R})$}{spin(1,7, R)} Triality}
In the case of an $8$ dimensional space-\textit{time} with symmetry group $\so(1,7,\mathbb{R})$, one of course will also find triality of representations. It should be noted most of the bases constructed here can be made in a quicker `ad-hoc' way by simply Wick rotating the right basis elements in both vector and spinor cases so the correct new metric signature is respected; for completeness a construction via Clifford algebras is still done for the spinor reps.\footnote{While we are keeping pedantry in the footnotes, notice the center of the connected component of $\grp{SO}(1,7)$ is trivial, and so in this case we do find spinor reps worthy of being prefix-less, as the constructed $\grp{Spin}(1,7)$ reps have a center of $\mathbb{Z}_2$.} We will also reuse the symbols $V, L, R$ for the representations. 

\subsection{Vector Representation}
The defining representation of $\grp{SO}(1,7)$ are those $8 \times 8$ real matrices $R$ satisfying:
\begin{equation}
    R^\intercal \eta R = \eta
\end{equation}
With $\eta = \mathbb{I}_{1,7}$ the flat metric of $8$ dimensional Minkowski spacetime (the $0^{\textrm{th}}$ axis is timelike). The standard basis for the corresponding representation of the groups Lie algebra may be given in a condensed form via
\begin{equation}
 \left(V_{ij}\right)^a{}_b = \eta_{j b} \delta^a{}_i - \eta_{i b} \delta^a{}_j
\end{equation}
with $\delta$ the Kronecker delta, and all indices run from $0$ to $7$ inclusive. The matrices $V_{ij}$ generate rotations/boosts. The span of these matrices along with the matrix commutator gives a representation of $\spin(1,7)$, which we will call the $\mathbf{8}_V$ representation. Finally we will alter this simple basis of rotation generators as before: 
\begin{equation}
    V_{1,5} \mapsto -V_{1,5}  \quad \mathrm{and} \quad  V_{2,6} \mapsto -V_{2,6},
\end{equation} 
all else left the same.

\subsection{Complex Chiral Spinor Representations}
Following from the basis for $C\ell(7,0)$ in Eq.~(\ref{eq:so(7)}), we may build a basis for the real Clifford algebra $C\ell(1,7)$ as follows:
\begin{equation}
\Gamma_0 = \begin{pmatrix}
\mathbb{I}_8 & 0 \\
0 & -\mathbb{I}_8
\end{pmatrix}, \quad \Gamma_{\mu} = i \begin{pmatrix}
0 & g_\mu \\
g_\mu & 0
\end{pmatrix} \quad 1 \leq \mu \leq 7.
\end{equation}
One can verify these satisfy the defining Clifford algebra relations: $\{\Gamma_i, \Gamma_j\} = 2 \eta_{ij}$. This representation is not manifestly chiral as in the compact case, and it is easy to compute the volume element $\omega = \prod^{8}_{i} \Gamma_i$ is a pseudo-scalar: $\omega ^2 = -1$, implying the chiral spinor representations are complex or quaternionic, and not real as in the compact case. To go to the chiral basis one  employs the following unitary the change of basis
\begin{equation}
    \Gamma_i \> \mapsto \>A \>\Gamma_i \> A^\dagger, \qquad A = \frac{1}{\sqrt{2}} \begin{pmatrix}
        \mathbb{I}_8 & -i \mathbb{I}_8 \\
        -i \mathbb{I}_8 & \mathbb{I}_8
    \end{pmatrix} .
\end{equation}
Once in the chiral basis the $\alg{spin}(1,7,\mathbb{R})$ algebra is generated as always by the second degree elements of the Clifford algebra basis (reusing notation):
\begin{equation}
\mathbf{8}_L \oplus \mathbf{8}_R = L^\prime_{ij} \oplus R^\prime_{ij}=\operatorname{span}\{\Gamma_i \Gamma_j/2\}, \quad 0 \leq i < j \leq 7.
\end{equation}
Taking the left and right pieces separately we will modify the bases in the following way:\footnote{The purpose of these modifications to the bases are not only to make the triality map simple, but also to ensure the constructed bases are in agreement with the previous constructions when considering identical sub-algebras, e.g. removing the $0^{th}$ axis in both cases.}
\begin{equation}
    L_{ij} \> = \> - M L^{\prime}_{ij} M^\dagger, \quad R_{ij} \> = \> - M^\dagger R^{\prime}_{ij} M,  \qquad M = \grp{diag}\{i,1,1,1,1,1,1,1\}.
\end{equation}
Finally, flip the sign again on those generators with $ij = \{15,26\}$. \textit{Now} it is true that the commutation relations of our Lorentzian vector generators, and of our two spinor representations, are respectively identical (i.e. our Lie algebra bases all have exactly the same structure constants), so we are in a most convenient position to build an algebra outer-automorphism. It is worth noting that for our spinor representations the generators of boosts are Hermitian, while rotations they are anti-Hermitian, just as for the vectors we have that boost generators are symmetric while rotations are anti-symmetric. All three of these representations have a simple Killing form, in the above bases given as:
\begin{equation}
    \kappa_{a b} \propto \begin{pmatrix}
        \mathbb{I}_{7} & 0 \\
        0 & -\mathbb{I}_{21}
    \end{pmatrix}.
\end{equation}
\subsubsection{Spinor duality: Complex Conjugation}
Before we look at the triality, we can ask about spinor duality: are these spinor representations quaternionic or complex? It is relatively straightforward from definitions to see in the bases given that the spinor representations are exactly complex conjugates of one another:
\begin{equation}
    L_{ij}^* = R_{ij} \quad \forall i,j.
\end{equation}
It is slightly more work, but simple algebra, to verify that the bases for $R$ and $L$ are not related by a change of basis, hence the representations are complex. As the vector representation constructed is manifestly real:
\begin{equation}
    V_{ij}^* = V_{ij} \quad \forall i,j,
\end{equation}
we thus have a $\mathbb{Z}_2$ action on the three $8$ dimensional representations, which swaps the spinor reps and acts trivially on the vector rep!

\section{Triality Map for \texorpdfstring{$\spin(1,7)$}{spin(1,7)}}
Triality acts similarly in the Lorentzian case, upon the same setup for quartets of generators defined in Eq.~(\ref{TrialityQuartets}), only the names now refer to the generators of $\spin(1,7)$. Define a matrix $T$ analogous to $H$ previously, as follows:
\begin{equation}
    T = \frac{1}{2}\begin{pmatrix}
        -1 & i & -i & -i \\
        i & 1 & 1 & 1 \\
        -i & 1 & 1 & -1 \\
        -i & 1 & -1 & 1 \\
    \end{pmatrix}.
\end{equation}
$T$ acting simultaneously on these quartets will send $V \rightarrow L \rightarrow R \rightarrow V$, not only at the level of the representation, but exactly to the same specified bases. We can see $T^3 = \mathbb{I}_4$. The reason $\spin(1,7)$ is not said to have triality exactly as $\spin(8)$ or $\spin(4,4)$ do, is because this triality map $T$ takes complex combinations of Lie algebra basis elements (as it must if $\mathbf{8}_V$ is real while $\mathbf{8}_L$, $\mathbf{8}_R$ are complex). In some sense as long as it is understood as a map between algebra representation basis elements and not an expression of vector addition, we can claim not to have broken any rules. Despite these complex combinations, because it is an algebra automorphism it also leaves the Killing form untouched, so we have not `accidentally complexified' at any point, and each representation is still an honest to goodness representation of $\spin$(1,7).

We can see the generators $T$ and $^*$ of our outer automorphism group map us directly from representation to representation, and conveniently precisely from specified basis to specified basis. It is quick to verify $T^2 = T^* = T^{-1}$, and so $T$ and complex conjugation do in fact generate the outer automorphism group $S_3$ of the Lie algebra. 
\subsection{Diagonal Triality: a \texorpdfstring{$\mathbb{Z}_3$}{Z3} grading of \texorpdfstring{$\spin(1,7)$}{spin(1,7)}}
Notice $T$ is symmetric. The normalized eigenvectors of $T$ are:
\begin{equation}
\label{Teigens}
    \begin{split}
    &\left| 1\right\rangle_3 = \frac{1}{\sqrt{2}}\begin{pmatrix}0 \\ 1 \\ 0 \\ 1\end{pmatrix}, \quad \left| 1\right\rangle_8 = \frac{1}{\sqrt{6}}\begin{pmatrix} 0 \\1 \\2 \\ -1\end{pmatrix}, \quad \left|e^{\pm i 2 \pi /3}\right\rangle = \frac{1}{\sqrt{6}}\begin{pmatrix} \pm \sqrt{3} \\ -1 \\ 1 \\ 1 \end{pmatrix}.
    \end{split}
\end{equation}
Recall that these eigenvectors are in the space of quartets of Lie algebra generators and as there are seven of these quartets, each of these eigenspaces is seven dimensional. Diagonalizing $T$ requires rewriting the bases of of the Lie algebra such that they differ from the easily interpretable coordinate plane of rotation basis. Still, each of these eigenvectors is a real vector, and so performing this change of Lie algebra basis will not alter the signature of the Killing form, and so we can diagonalize the action of $T$ without leaving $\spin(1,7)$. This gives a real-orthogonal change of basis on the $28$ dimensional Lie algebra as $$B = \left(\left| 1\right\rangle_3 \> \left| 1\right\rangle_8 \> \left| e^{i 2 \pi /3}\right\rangle \> \left| e^{-i 2 \pi /3}\right\rangle\right)$$ acting on the quartets. That is to say we construct a new basis of generators for $V$, $L$, and $R$ given by the image of $B$ on the quartets of generators. 

After diagonalizing, $T = \operatorname{diag}\{1,1,e^{i 2 \pi/3},e^{-i 2 \pi/3}\}$. The $28$ dimensional algebra of $\spin(1,7)$ is broken into three pieces by the eigenspaces under triality: the fourteen dimensional triality invariant subspace can be seen to be the compact sub-algebra $\alg{g}_2$, given in an $8$ dimensional real representation that looks like $1\oplus7$. The other 14 vectors in the algebra form `conjugate pairs' of generators, related by an anti-linear automorphism $\theta$ given in the next section. These 14 generators are each null with respect to the Killing form, and orthogonal to all other vectors aside from their $\theta$ conjugate. The Killing form in this basis looks as:
\begin{equation}
    \kappa_{a b} \propto \begin{pmatrix}
        -\mathbb{I}_{14} & 0_{7} & 0_{7} \\
        0_{7} & 0_{7} & -\mathbb{I}_7 \\
        0_7 & -\mathbb{I}_7 & 0_7
    \end{pmatrix}.
\end{equation}
The algebra may be decomposed into the first 14 bases, the next 7 bases $N^+$, and the last 7 bases $N^-$:
\begin{equation}
    \spin(1,7) \equiv \left| e^{-i 2\pi/3}\right\rangle \oplus \left| 1\right\rangle \oplus \left| e^{+i 2\pi/3}\right\rangle = N^- \oplus \alg{g}_2 \oplus N^+
\end{equation}
Complex conjugation of course will send $N^{\pm}_L \leftrightarrow N^{\pm}_R$ while fixing the vector rep. As for triality, given a particular representation such as $\mathbf{8}_L$, triality acts on this decomposition diagonally as:
\begin{equation}
    \mathbf{8}_R =T\cdot \mathbf{8}_L \equiv \left(e^{-i2\pi/3} N_L^-\right) \oplus \alg{g}_2 \oplus \left(e^{i2\pi/3} N_L^+ \right) = N_R^- \oplus \alg{g}_2 \oplus N_R^+,
\end{equation}
I.e. the only difference between the three representations is a phase upon the null-vector generators. This of course \textit{only} works with null generators because if the generator were not null the signature of the Killing form would change after application of a complex phase. 

The triality eigenvalue characterizes the commutation relations of this decomposition, since the phase is multiplicative as well as an algebra automorphism: the commutator of two elements with triality eigenvalues $\lambda_1$, $\lambda_2$ gives an element with a corresponding eigenvalue given by the product $\lambda_1 \cdot \lambda_2$. From this it is easy to see:
\begin{equation}
\begin{split}
    \left[\alg{g}_2,\alg{g}_2\right] = \alg{g}_2 , &\quad \left[\alg{g}_2,N^{\pm}\right] = N^{\pm}, \\
    \left[N^+,N^-\right] = \alg{g}_2, &\quad \left[N^{\pm},N^{\pm}\right] = N^{\mp}, 
\end{split}
\end{equation}
It is this last commutator being non-zero which distinguishes this algebra as being $\mathbb{Z}_3$ graded, as opposed to merely three graded. 
\subsubsection{Comment on Jordan Algebra Structure}
In this decomposition we have an anti-linear relation $\theta$ between $N^+$ and $N^-$, which is valid in all three representations, though the anti-linearity is not manifest in the vector rep:
\begin{equation}
    N^- = \theta(N^+) = \tau^{-1} (N^+)^* \tau, \quad \mathrm{with}\quad \tau = \operatorname{diag}\{-1,1,1,1,1,1,1,1\} .
\end{equation}
Where $*$ means entry-wise complex conjugation on a given basis. As each representation preserves the bilinear inner product defined by $\eta$, which has up to a sign the same components of $T$, these relations may also be phrased another way: 
\begin{equation}
    N^- =  - (N^+)^\dagger.
\end{equation}
The action of $\theta$ can be extended to the whole algebra because $\theta$ acts as the identity of the given real representation of the $\alg{g}_2$ subalgebra. If we instead consider $\theta$ as an active transformation of the entire algebra, it will send the left and right handed representations to one another, followed by a time reversal change of basis:
\begin{equation}
    N^-_L = \tau^{-1}(N^+_L)^* \tau =  \tau^{-1}(N^+_R) \tau.
\end{equation}
With the action of $\theta$ understood, the $\mathbb{Z}_3$-grading gives a Jordan algebra\cite{Chu_2011} upon $N^+$:
\begin{equation}
    \left[\left[x,\theta(y)\right],z\right], \quad \forall x,y,z\in N^+
\end{equation}
However in this case we have other options for the Jordan product, such as: 
\begin{equation}
\label{twistedJ}
    \left[\left[x,y\right],\theta(z)\right], \quad \forall x,y,z\in N^+,
\end{equation}
twisting any component works just as well in this case.  If the Killing form is twisted sesqui-linearly to the form $\operatorname{Tr}(x^\dagger y)$ for $x,y\in N^+$ this can serve as a Hermitian inner product on $N^+$. It is my suspicion that the existence of these Jordan products on the $N$ spaces is merely the statement that there is a highly non-unique cross product in $7$ dimensions. But perhaps there is more to say. I believe this Jordan algebra has no tripotents, and of course for Eq.~(\ref{twistedJ}) there can be none. 

\section{Final Remarks}
In terms of the group and representation theory, $\so(8)$ is very well understood, but as for applications many open questions still remain. Firstly, whether there is anything interesting in the aforementioned possible triality of twistors in $1+7$ dimensional spacetime. It is interesting that here we find a real vector representation where triality acts diagonally on the basis of generators, and so in some sense the transformations of spinors in $8$ dimensional spacetime are merely a phase (\textit{something} like a Wick rotation) applied to our fourteen null generators. I.e. the only difference between the same transformation from $N^+$ applied to a vector $v$ and a spinor $s$ of $\grp{Spin}(1,7)$:
\begin{equation}
    e^{V_{ij} \theta_{ij}}v, \quad \mathrm{vs.} \quad  e^{L_{ij} \phi_{ij}}v
\end{equation}
is that $\phi_{ij} = e^{i 2\pi/3}\theta_{ij}$. So somehow `Wick' rotating not by $\pi/2$ but by $2\pi/3$ can map, in certain cases, between a vector and spinor. I'm not really certain this makes sense as a question to investigate further, but it feels like it \textit{almost} does. With $\grp{Spin}(1,7)$ as the conformal group of $6$ dimensional space, perhaps there is some relation between the conformal light-cone here, this seven dimensional space of null vectors, triality, and $SU(4) \hookleftarrow U(1)\times SU(3)$? There is likely \textit{a} story here, but I am not sure what it has to do with reality. 

Secondly, the most tantalizing open question about triality is of course whether it can be used to explain the three generations of matter: the duality of field theory solutions due to $CPT$ invariance leads to a doubling of particle content in the standard model, namely anti-matter. It has been considered well before myself, whether the three fold symmetry of $\so(8)$ might be able explain the three fold copies of matter we see in nature, the generations. One obstruction to making this work in the `obvious' way is finding outer automorphisms which act at the level of the representation, and not merely on the algebra. For example in the Euclidean case $K$ and $H$ both act on the algebra, and so the automorphism has no `intertwining' action upon the representation space (the fields). However in the Lorentzian case with $T$ and $^*$, complex conjugation \textit{can} be made to act not only upon the Lie algebra generators, but upon on the fields as well. Thus the question becomes if there is there some way, some combination of automorphisms, in some base field for our representation or some hyper-complex number system, where we can have a non-linear intertwiner which cubes to the identity, yet its square remains a non-linear intertwiner as well? For example a standard conjugate-linear operator cannot work, because when squared it becomes linear, and so this could not garner us three distinct representations. Perhaps if one can realizes the triality as three dualities it could be feasible.

This work is foremost a pedagogical piece --- it is the very kind of paper I wish existed when I first set out to understand triality, to get my hands on it, as it were, and see for myself what `performing the triality' actually looks like. I hope new students and experts alike can find some useful insights from what I have learned in my jaunt through the subject, and if nothing else, be brought some convenience through the existence of the constructions in the paper. 

\newpage
\bibliographystyle{plain} % We choose the "plain" reference style
\bibliography{TrialityMcRae.bib} % Entries are in the refs.bib file

\end{changemargin}
\end{document}